\documentclass[10pt,letterpaper]{article}
\usepackage[top=0.85in,left=1in,footskip=0.75in,marginparwidth=2in]{geometry}

\usepackage[utf8]{inputenc}
\usepackage{amsmath}
\usepackage{cite}

\usepackage{nameref,hyperref}

\usepackage[right]{lineno}

\usepackage{microtype}
\DisableLigatures[f]{encoding = *, family = * }

\usepackage{changepage}

\usepackage[aboveskip=1pt,labelfont=bf,labelsep=period,singlelinecheck=off]{caption}

\makeatletter
\renewcommand{\@biblabel}[1]{\quad#1.}
\makeatother

\usepackage{lastpage,fancyhdr,graphicx}
\usepackage{epstopdf}
\fancyheadoffset[L]{2.25in}
\fancyfootoffset[L]{2.25in}

\usepackage{color}

\definecolor{Gray}{gray}{.25}

\usepackage{graphicx}

\usepackage{sidecap}

\usepackage{wrapfig}
\usepackage[pscoord]{eso-pic}
\usepackage[fulladjust]{marginnote}
\reversemarginpar

\begin{document}
\vspace*{0.35in}

\begin{flushleft}
{\Large
\textbf{Addition formulas for the $\boldsymbol{_pF_p}$ and $\boldsymbol{_{p+1}F_p}$ generalized hypergeometric functions with arbitrary parameters and their Kummer- and Euler-type transformations}
}
\newline
\\
Krishna Choudhary
\\
\bigskip
Gladstone Institutes, 1650 Owens St, San Francisco, CA 94158
\\
\bigskip
krishna.choudhary@gladstone.ucsf.edu, kchoudhary@ucdavis.edu

\end{flushleft}

\section*{Abstract}
We obtain addition formulas for $_pF_p$ and $_{p+1}F_p$ generalized hypergeometric functions with general parameters. These are utilized in conjunction with integral representations of these functions to derive Kummer- and Euler-type transformations that express $_pF_p\left(x\right)$ and $_{p+1}F_p\left(x\right)$ in the form of sums of $_pF_p\left(-x\right)$ and $_{p+1}F_p\left(-x\right)$ functions, respectively.

\linenumbers

\section{Introduction}

A generalized hypergeometric function with $p$ numerator parameters and $q$ denominator parameters is defined as 
\begin{eqnarray}
_pF_q\left( 
\begin{smallmatrix}
a_1, & a_2, & \ldots, & a_p\\
b_1, & b_2, & \ldots, & b_q \\
\end{smallmatrix};x
\right) & = & \sum_{i=0}^\infty \frac{\left(a_1\right)_i \left(a_2\right)_i \ldots \left(a_p\right)_i}{\left(b_1\right)_i \left(b_2\right)_i \ldots \left(b_q\right)_i}\frac{x^i}{i!},
\end{eqnarray}
where the parameters and the argument can take complex values except that the denominator parameters cannot be negative integers, and $\left(a\right)_i$ is the Pochhammer symbol for the ascending factorial $\prod_{j=0}^{i-1}\left(a+j\right)$ with $\left(a\right)_0 = 1$. For example, confluent hypergeometric function with $p=q=1$ and Gaussian hypergeometric function with $p=q+1=2$ are two of its special cases that are well-known. 

The generalized hypergeometric functions are of significant interest for their applications in diverse areas, notably including mathematical physics and mathematical statistics \cite{slater1966generalized, seaborn2013hypergeometric}. Insights into their theory as well as their applications are facilitated by various transformation formulas and identities, many of which have been collected in standard references \cite{slater1966generalized, bateman1953higher, abramowitz1970handbook, beals2016special}. As mathematical models in physics or statistics grow in complexity and involve systems that are described in terms of multiple parameters, their numerical or analytical solutions may have to contend with cases where $p,q\geq2$. In the last two decades, a number of transformation and summation formulas have been derived for such cases. Paris utilized the addition theorem for $_1F_1\left(x\right)$ in conjunction with an integral representation of hypergeometric functions to derive a Kummer-type transformation formula that connects $_2F_2\left(x\right)$ with general parameters to $_2F_2\left(-x\right)$ \cite{paris2005kummer}. A number of formulas for functions with special relationships between the numerator and denominator parameters or special values of the argument have also been derived \cite{miller2011euler, kim2012two, rakha2014extension, srivastava2019extensions}. For example, Miller and Paris derived Kummer- and Euler-type transformation formulas respectively for $_pF_p\left(x\right)$ and $_{p+1}F_p\left(x\right)$ with integral differences between the numerator and denominator parameters \cite{miller2013transformation}. However, to the best of my knowledge, such formulas are not available for $_pF_p\left(x\right)$ and $_{p+1}F_p\left(x\right)$ with general parameters. If addition formulas were available for the $_pF_p$ and $_{p+1}F_p$  functions, one could follow the method of Paris \cite{paris2005kummer} and derive the said transformation formulas. However, to the best of my knowledge, addition formulas, which are important in their own right \cite{koelinkEOM}, are also not available for the $_pF_p$ and $_{p+1}F_p$ functions. Here, I fill this gap by stating and proving addition formulas for these functions, which I use in conjunction with their integral representations to derive transformation formulas that connect $_pF_p\left(x\right)$ and $_{p+1}F_p\left(x\right)$ with $_pF_p\left(-x\right)$ and $_{p+1}F_p\left(-x\right)$, respectively for general parameter values. 

\section{Theorem 1 (addition formula for $\boldsymbol{_pF_p}$).} For $p \geq 1$ and $b_i \neq 0, -1, -2, ...$ $\forall$ $1 \leq i \leq p$, 
\begin{eqnarray}
_{p}F_{p}\left( 
\begin{smallmatrix}
a_1 , & a_2, & ..., & a_p \\
b_1 , & b_2, & ..., & b_p  \\
\end{smallmatrix};x+y
\right) & = & e^x \sum_{j_1=0}^\infty \sum_{j_2=0}^\infty ... \sum_{j_p=0}^\infty \prod_{q=1}^p \frac{\left(b_q - a_q\right)_{j_q} \left(a_q\right)_{u_{q-1}}}{\left(b_q\right)_{u_q}} \frac{\left(-x\right)^{j_q}}{j_q!} \cdot \nonumber \\
  & & \hspace{35pt} {_{p}}F_{p}\left( 
\begin{smallmatrix}
a_1 + u_0, & a_2 + u_1, & a_3 + u_2, & ..., & a_p + u_{p-1}\\
b_1 + u_1, & b_2 + u_2, & b_3 + u_3, & ..., & b_p + u_p \\
\end{smallmatrix};y
\right), \label{eq: addition_theorem_pFp}
\end{eqnarray}
where $u_0 = 0$, $u_1 = j_1$, $u_2 = j_1+j_2$, ..., $u_p = \sum_{q=1}^p j_q$.

\subsubsection*{Proof} We can use the method of induction to prove this theorem. First, note that for $p=1$, Eq. \ref{eq: addition_theorem_pFp} becomes 
\begin{eqnarray}
_{1}F_{1}\left( 
\begin{smallmatrix}
a_1 \\
b_1 \\
\end{smallmatrix};x+y
\right) & = & e^x \sum_{j_1=0}^\infty \frac{\left(b_1 - a_1\right)_{j_1}}{\left(b_1\right)_{j_1}} \frac{\left(-x\right)^{j_1}}{j_1!} \cdot {_{1}}F_{1}\left( 
\begin{smallmatrix}
a_1 \\
b_1 + j_1 \\
\end{smallmatrix};y
\right), 
\end{eqnarray}
which is correct (see Eq. 2.3.5 in Slater \cite{slater1960confluent}). Next, let us assume that Eq. \ref{eq: addition_theorem_pFp} is valid for some $p=k$. Then, using the integral representation in Eq. 4.8.3.12 of Slater \cite{slater1966generalized}, for $p=k+1$ and $\text{Re}\left(b_{k+1}\right) > \text{Re}\left(a_{k+1}\right)>0$,
\begin{align}
& _{k+1}F_{k+1}\left( 
\begin{smallmatrix}
a_1, & ..., & a_{k+1} \\
b_1, & ..., & b_{k+1}  \\
\end{smallmatrix};x+y
\right) \nonumber \\
& = \frac{\Gamma\left(b_{k+1}\right)}{\Gamma\left(a_{k+1}\right) \Gamma\left(b_{k+1} - a_{k+1}\right)} \int_0^1 t^{a_{k+1} - 1} \left(1-t\right)^{b_{k+1} - a_{k+1} - 1} {_{k}}F_{k}\left( 
\begin{smallmatrix}
a_1, & ..., & a_{k} \\
b_1, & ..., & b_{k}  \\
\end{smallmatrix};xt + yt
\right) dt. \label{eq: pFp_add_thrm_proof_start}
\end{align}
Here, using Eq. \ref{eq: addition_theorem_pFp} for $p=k$ given $b_i \neq 0, -1, -2, ...$ $\forall$ $1 \leq i \leq k$ yields
\begin{align}
& \frac{\Gamma\left(b_{k+1}\right)}{\Gamma\left(a_{k+1}\right) \Gamma\left(b_{k+1} - a_{k+1}\right)} \int_0^1 t^{a_{k+1} - 1} \left(1-t\right)^{b_{k+1} - a_{k+1} - 1} \cdot \nonumber \\
& \hspace{35pt} e^{xt} \sum_{j_1=0}^\infty \sum_{j_2=0}^\infty ... \sum_{j_k=0}^\infty \prod_{q=1}^k \frac{\left(b_q - a_q\right)_{j_q} \left(a_q\right)_{u_{q-1}}}{\left(b_q\right)_{u_q}} \frac{\left(-xt\right)^{j_q}}{j_q!} {_{k}}F_{k}\left( 
\begin{smallmatrix}
a_1 + u_0, & ..., & a_k + u_{k-1}\\
b_1 + u_1, & ..., & b_k + u_k \\
\end{smallmatrix};yt
\right) dt, \nonumber
\end{align}
where $u_0 = 0$, $u_1 = j_1$, $u_2 = j_1+j_2$, ..., $u_k = \sum_{q=1}^k j_q$. Upon rearrangement of terms and reversal of the order of summation and integration, this expression can be written as
\begin{align}
& e^x \frac{\Gamma\left(b_{k+1}\right)}{\Gamma\left(a_{k+1}\right) \Gamma\left(b_{k+1} - a_{k+1}\right)} \sum_{j_1=0}^\infty \sum_{j_2=0}^\infty ... \sum_{j_k=0}^\infty \prod_{q=1}^k \frac{\left(b_q - a_q\right)_{j_q} \left(a_q\right)_{u_{q-1}}}{\left(b_q\right)_{u_q}} \frac{\left(-x\right)^{j_q}}{j_q!} \cdot \nonumber \\
& \hspace{35pt} \int_0^1 t^{a_{k+1} + u_k - 1} \left(1-t\right)^{b_{k+1} - a_{k+1} - 1} e^{-x\left(1-t\right)} {_{k}}F_{k}\left( 
\begin{smallmatrix}
a_1 + u_0, & ..., & a_k + u_{k-1}\\
b_1 + u_1, & ..., & b_k + u_k \\
\end{smallmatrix};yt
\right) dt. \nonumber
\end{align}
Next, using the power series expansion for $e^{-x\left(1-t\right)}$ in the above expression, we write it as
\begin{align}
& e^x \frac{\Gamma\left(b_{k+1}\right)}{\Gamma\left(a_{k+1}\right) \Gamma\left(b_{k+1} - a_{k+1}\right)} \sum_{j_1=0}^\infty \sum_{j_2=0}^\infty ... \sum_{j_k=0}^\infty \prod_{q=1}^k \frac{\left(b_q - a_q\right)_{j_q} \left(a_q\right)_{u_{q-1}}}{\left(b_q\right)_{u_q}} \frac{\left(-x\right)^{j_q}}{j_q!} \cdot \nonumber \\
& \int_0^1 t^{a_{k+1} + u_k - 1} \left(1-t\right)^{b_{k+1} - a_{k+1} - 1} \sum_{j_{k+1}=0}^\infty \frac{\left(-x\right)^{j_{k+1}} \left(1-t\right)^{j_{k+1}}}{j_{k+1}!} {_{k}}F_{k}\left(
\begin{smallmatrix}
a_1 + u_0, & ..., & a_k + u_{k-1}\\
b_1 + u_1, & ..., & b_k + u_k \\
\end{smallmatrix};yt
\right) dt, \nonumber
\end{align}
which upon reversal of the order of integration and summation yields
\begin{align}
& e^x \frac{\Gamma\left(b_{k+1}\right)}{\Gamma\left(a_{k+1}\right) \Gamma\left(b_{k+1} - a_{k+1}\right)} \sum_{j_1=0}^\infty \sum_{j_2=0}^\infty ... \sum_{j_k=0}^\infty \prod_{q=1}^k \frac{\left(b_q - a_q\right)_{j_q} \left(a_q\right)_{u_{q-1}}}{\left(b_q\right)_{u_q}} \frac{\left(-x\right)^{j_q}}{j_q!} \cdot \nonumber \\
& \sum_{j_{k+1}=0}^\infty \frac{\left(-x\right)^{j_{k+1}} }{j_{k+1}!} \int_0^1 t^{a_{k+1} + u_k - 1} \left(1-t\right)^{b_{k+1} - a_{k+1} +j_{k+1} - 1}  {_{k}}F_{k}\left(
\begin{smallmatrix}
a_1 + u_0, & ..., & a_k + u_{k-1}\\
b_1 + u_1, & ..., & b_k + u_k \\
\end{smallmatrix};yt
\right) dt. \nonumber
\end{align}
Once again, by using the integral representation in Eq. 4.8.3.12 of Slater \cite{slater1966generalized} in the above expression, we get
\begin{align}
& e^x \frac{\Gamma\left(b_{k+1}\right)}{\Gamma\left(a_{k+1}\right) \Gamma\left(b_{k+1} - a_{k+1}\right)} \sum_{j_1=0}^\infty \sum_{j_2=0}^\infty ... \sum_{j_k=0}^\infty \prod_{q=1}^k \frac{\left(b_q - a_q\right)_{j_q} \left(a_q\right)_{u_{q-1}}}{\left(b_q\right)_{u_q}} \frac{\left(-x\right)^{j_q}}{j_q!} \cdot \nonumber \\
& \sum_{j_{k+1}=0}^\infty \frac{\left(-x\right)^{j_{k+1}} }{j_{k+1}!} \frac{\Gamma\left(b_{k+1} - a_{k+1} + j_{k+1}\right) \Gamma\left(a_{k+1} + u_k\right)}{\Gamma\left(b_{k+1} + u_{k+1}\right)} \cdot {_{k+1}}F_{k+1}\left( 
\begin{smallmatrix}
a_1 + u_0, & ..., & a_{k+1} + u_k\\
b_1 + u_1, & ..., & b_{k+1} + u_{k+1}, \\
\end{smallmatrix};y
\right), \nonumber
\end{align}
where $u_{k+1} = \sum_{q=1}^{k+1}j_{q}$. This expression can be more compactly written as
\begin{align}
& e^x \sum_{j_1=0}^\infty \sum_{j_2=0}^\infty ... \sum_{j_{k+1}=0}^\infty \prod_{q=1}^{k+1} \frac{\left(b_q - a_q\right)_{j_q} \left(a_q\right)_{u_{q-1}}}{\left(b_q\right)_{u_q}} \frac{\left(-x\right)^{j_q}}{j_q!} \cdot {_{k+1}}F_{k+1}\left( 
\begin{smallmatrix}
a_1 + u_0, & ..., & a_{k+1} + u_{k}\\
b_1 + u_1, & ..., & b_{k+1} + u_{k+1} \\
\end{smallmatrix};y
\right) \nonumber
\end{align}
to get the right hand side of Eq. \ref{eq: addition_theorem_pFp} for $p=k+1$. While we have shown that the left hand side in Eq. \ref{eq: pFp_add_thrm_proof_start} is equal to the above expression given that $\text{Re}\left(b_{k+1}\right) > \text{Re}\left(a_{k+1}\right)>0$, they are equal for general values of $b_{k+1}$ and $a_{k+1}$ as well due to analytic continuation since they are both analytic functions of these parameters, given $b_{k+1} \neq 0, -1, -2, ...$. This implies that if the theorem holds for $p=k$, it also holds for $p=k+1$. Since we know that it holds for $p=1$, by induction, it holds for all $p \geq 1$ given $b_{i} \neq 0, -1, -2, ...$ $\forall$ $1 \leq i \leq p$, thereby completing the proof.

\section{Theorem 2 (addition formula for $\boldsymbol{_{p+1}F_p}$).} For $p \geq 1$ and given $b_i \neq 0, -1, -2, ...$ $\forall$ $1 \leq i \leq p$, $\left|x+y\right| < 1$, $\left|y\right| < \left|x\right|$ and $\text{Re}\left(x\right) < 1/2$,
\begin{align}
& _{p+1}F_{p}\left( 
\begin{smallmatrix}
a_0, & a_1 , & ..., & a_p \\
& b_1 , & ..., & b_p  \\
\end{smallmatrix};x+y
\right) \nonumber \\
& = \left(\frac{1}{1-x}\right)^{a_0} \sum_{j_1 = 0}^\infty \sum_{j_2 = 0}^\infty ... \sum_{j_p = 0}^\infty \left(a_0\right)_{u_p} \prod_{q=1}^p \frac{\left(b_q - a_q\right)_{j_q} \left(a_q\right)_{u_{q-1}}}{\left(b_q\right)_{u_q} j_q !} \left(\frac{x}{x-1}\right)^{j_q}  \cdot \nonumber \\
  & \hspace{25pt} {_{p+1}}F_{p}\left(
\begin{smallmatrix}
a_0+u_p, & a_1 + u_0, & a_2 + u_1, & a_3 + u_2, & ..., & a_p + u_{p-1}\\
& b_1 + u_1, & b_2 + u_2, & b_3 + u_3, & ..., & b_p + u_p\\
\end{smallmatrix};\frac{-y}{x-1}
\right), \label{eq: add_theorem_p+1Fp}
\end{align}
where $u_0 = 0$, $u_1 = j_1$, $u_2 = j_1+j_2$, ..., $u_p = \sum_{q=1}^p j_q$.

\subsubsection*{Proof} The addition formula for $_{p+1}F_{p}$ can be derived using the addition formula for $_{p}F_{p}$ by appealing to a relationship between the two derived by application of Mellin transform (see Eq. 4.8.3.3 in Slater \cite{slater1966generalized}),
\begin{eqnarray}
_{p+1}F_{p}\left( 
\begin{smallmatrix}
a_0, & a_1 , & a_2 , & ..., & a_p \\
& b_1 , & b_2 , & ..., & b_p  \\
\end{smallmatrix};x+y
\right) & = & \frac{1}{\Gamma\left(a_0\right)} \int_0^\infty t^{a_0 -1} e^{-t} {_p}F_{p}\left( 
\begin{smallmatrix}
a_1 , & a_2, & ..., & a_p \\
b_1 , & b_2, & ..., & b_p  \\
\end{smallmatrix};xt+yt
\right) dt \nonumber
\end{eqnarray}
given $\left|x+y\right| < 1$. Now, we work with the right hand side of the above equation and show that it equals the right hand side of Eq. \ref{eq: add_theorem_p+1Fp}. Given $b_i \neq 0, -1, -2, ...$ $\forall$ $1 \leq i \leq p$, we use Eq. \ref{eq: addition_theorem_pFp} to rewrite the right hand side of the above equation as
\begin{align}
& \frac{1}{\Gamma\left(a_0\right)} \int_0^\infty t^{a_0 -1} e^{-t} e^{xt} \sum_{j_1=0}^\infty \sum_{j_2=0}^\infty ... \sum_{j_p=0}^\infty \prod_{q=1}^p \frac{\left(b_q - a_q\right)_{j_q} \left(a_q\right)_{u_{q-1}}}{\left(b_q\right)_{u_q}} \frac{\left(-xt\right)^{j_q}}{j_q!} \cdot \nonumber \\
  & \hspace{55pt} {_{p}}F_{p}\left( 
\begin{smallmatrix}
a_1 + u_0, & a_2 + u_1, & ..., & a_p + u_{p-1}\\
b_1 + u_1, & b_2 + u_2, & ..., & b_p + u_p \\
\end{smallmatrix};yt
\right) dt, \nonumber 
\end{align}
where $u_0 = 0$, $u_1 = j_1$, $u_2 = j_1+j_2$, ..., $u_p = \sum_{q=1}^p j_q$. Upon changing the order of integration and summation and using the series expansion for $_pF_p$, it yields
\begin{align}
& \sum_{j_1=0}^\infty \sum_{j_2=0}^\infty ... \sum_{j_p=0}^\infty \prod_{q=1}^p \frac{\left(b_q - a_q\right)_{j_q} \left(a_q\right)_{u_{q-1}}}{\left(b_q\right)_{u_q}} \frac{\left(-x\right)^{j_q}}{j_q!} \frac{1}{\Gamma\left(a_0\right)} \cdot \nonumber \\
& \hspace{35pt} \int_0^\infty t^{a_0 + u_p -1} e^{-t\left(1-x\right)} \sum_{i=0}^\infty \prod_{r=1}^{p} \frac{\left(a_r + u_{r-1}\right)_i}{\left(b_r + u_{r}\right)_i} \frac{y^i t^i}{i!} dt. \nonumber 
\end{align}
Once again, rearranging the terms and changing the order of integration and summation, we get
\begin{align}
& \sum_{j_1=0}^\infty \sum_{j_2=0}^\infty ... \sum_{j_p=0}^\infty \prod_{q=1}^p \frac{\left(b_q - a_q\right)_{j_q} \left(a_q\right)_{u_{q-1}}}{\left(b_q\right)_{u_q}} \frac{\left(-x\right)^{j_q}}{j_q!} \frac{1}{\Gamma\left(a_0\right)} \cdot \nonumber \\
& \hspace{35pt} \sum_{i=0}^\infty \prod_{r=1}^{p} \frac{\left(a_r + u_{r-1}\right)_i}{\left(b_r + u_{r}\right)_i} \frac{y^i}{i!} \int_0^\infty t^{a_0 + u_p + i -1} e^{-t\left(1-x\right)}  dt. \nonumber
\end{align}
Next, we substitute $v=t\left(1-x\right)$ in the integral to get
\begin{align}
& \sum_{j_1=0}^\infty \sum_{j_2=0}^\infty ... \sum_{j_p=0}^\infty \prod_{q=1}^p \frac{\left(b_q - a_q\right)_{j_q} \left(a_q\right)_{u_{q-1}}}{\left(b_q\right)_{u_q}} \frac{\left(-x\right)^{j_q}}{j_q!} \frac{1}{\Gamma\left(a_0\right)} \cdot \nonumber \\
& \hspace{35pt} \sum_{i=0}^\infty \prod_{r=1}^{p} \frac{\left(a_r + u_{r-1}\right)_i}{\left(b_r + u_{r}\right)_i} \frac{y^i}{i!} \left(\frac{1}{1-x}\right)^{a_0+u_p+i} \int_0^\infty v^{a_0 + u_p + i -1} e^{-v}  dv. \nonumber 
\end{align}
The integral in the above expression can be written as a gamma function yielding
\begin{align}
& \sum_{j_1=0}^\infty \sum_{j_2=0}^\infty ... \sum_{j_p=0}^\infty \prod_{q=1}^p \frac{\left(b_q - a_q\right)_{j_q} \left(a_q\right)_{u_{q-1}}}{\left(b_q\right)_{u_q}} \frac{\left(-x\right)^{j_q}}{j_q!} \frac{1}{\Gamma\left(a_0\right)} \cdot \nonumber \\
& \hspace{35pt} \sum_{i=0}^\infty \prod_{r=1}^{p} \frac{\left(a_r + u_{r-1}\right)_i}{\left(b_r + u_{r}\right)_i} \frac{y^i}{i!} \left(\frac{1}{1-x}\right)^{a_0+u_p+i} \cdot \Gamma\left(a_0 + u_p + i\right), \nonumber
\end{align}
which upon rearranging the terms and using $\left(a_0\right)_{u_p + i} = \left(a_0\right)_{u_p} \left(a_0 + u_p\right)_{i}$ can be written as
\begin{align}
& \left(\frac{1}{1-x}\right)^{a_0} \sum_{j_1=0}^\infty \sum_{j_2=0}^\infty ... \sum_{j_p=0}^\infty \left(a_0\right)_{u_p} \prod_{q=1}^p \frac{\left(b_q - a_q\right)_{j_q} \left(a_q\right)_{u_{q-1}}}{\left(b_q\right)_{u_q}} \left(\frac{x}{x-1}\right)^{j_q} \frac{1}{j_q!}  \cdot \nonumber \\
& \hspace{35pt} \sum_{i=0}^\infty \left(a_0+u_p\right)_{i} \prod_{r=1}^{p} \frac{\left(a_r + u_{r-1}\right)_i}{\left(b_r + u_{r}\right)_i} \frac{1}{i!} \left(\frac{-y}{x-1}\right)^{i}.  \nonumber
\end{align}
Finally, since $\left|y\right| < \left|x-1\right|$ given $\left|y\right| < \left|x\right|$ and $\text{Re}\left(x\right) < 1/2$, it can be expressed more compactly as
\begin{align}
& \left(\frac{1}{1-x}\right)^{a_0} \sum_{j_1 = 0}^\infty \sum_{j_2 = 0}^\infty ... \sum_{j_p = 0}^\infty \left(a_0\right)_{u_p} \prod_{q=1}^p \frac{\left(b_q - a_q\right)_{j_q} \left(a_q\right)_{u_{q-1}}}{\left(b_q\right)_{u_q} j_q !} \left(\frac{x}{x-1}\right)^{j_q}  \cdot \nonumber \\
  &  \hspace{25pt} {_{p+1}}F_{p}\left( 
\begin{smallmatrix}
a_0+u_p, & a_1 + u_0, & a_2 + u_1, & a_3 + u_2, & ..., & a_p + u_{p-1}\\
& b_1 + u_1, & b_2 + u_2, & b_3 + u_3, & ..., & b_p + u_p\\
\end{smallmatrix};\frac{-y}{x-1}
\right), \nonumber
\end{align}
which is the same as the right hand side of Eq. \ref{eq: add_theorem_p+1Fp}, thereby proving the theorem.

\section{Theorem 3 (Kummer-type transformation for $\boldsymbol{_pF_p}$ with general parameters).} For $p \geq 1$ and $b_i \neq 0, -1, -2, ...$ $\forall$ $1 \leq i \leq p+1$,
\begin{align}
& _{p+1}F_{p+1}\left(
\begin{smallmatrix}
a_1 , & a_2 , & a_3 , & ..., & a_p , & a_{p+1}\\
b_1 , & b_2 , & b_3 , & ..., & b_p , & b_{p+1} \\
\end{smallmatrix};x
\right) \nonumber \\
& = e^x \sum_{j_1=0}^\infty \sum_{j_2=0}^\infty ... \sum_{j_p=0}^\infty \prod_{q=1}^p \frac{\left(b_q - a_q\right)_{j_q} \left(a_q\right)_{u_{q-1}}}{\left(b_q\right)_{u_q}} \frac{\left(-x\right)^{j_q}}{j_q!} \cdot \nonumber \\
  & \hspace{55pt} {_{p+1}}F_{p+1}\left(
\begin{smallmatrix}
a_1 + u_0, & a_2 + u_1, & a_3 + u_2, & ..., & a_p + u_{p-1}, & b_{p+1} - a_{p+1}\\
b_1 + u_1, & b_2 + u_2, & b_3 + u_3, & ..., & b_p + u_p, & b_{p+1} \\
\end{smallmatrix};-x
\right), \label{eq: kummer_transform_pFp}
\end{align}
where $u_0 = 0$, $u_1 = j_1$, $u_2 = j_1+j_2$, ..., $u_p = \sum_{q=1}^p j_q$.

\subsubsection*{Proof} 
To prove the theorem, we utilize the approach followed by Paris \cite{paris2005kummer} for a Kummer-type transformation of $_2F_2\left(x\right)$. Hence, we utilize the integral representation for $_{p+1}F_{p+1}\left(x\right)$ function, which is available from Slater \cite{slater1966generalized} (see Eq. 4.8.3.12) and the addition formula for $_pF_p\left(x\right)$, which we obtained above (theorem 1). The integral representation for $_{p+1}F_{p+1}\left(x\right)$ is 
\begin{eqnarray}
_{p+1}F_{p+1}\left( 
\begin{smallmatrix}
a_1, & ..., & a_{p+1}\\
b_1, & ..., & b_{p+1}\\
\end{smallmatrix};x
\right) & = & \frac{\Gamma\left(b_{p+1}\right)}{\Gamma\left(a_{p+1}\right) \Gamma\left(b_{p+1}-a_{p+1}\right)} \nonumber \\
& & \hspace{-40pt} \int_0^1 t^{a_{p+1}-1} \left(1-t\right)^{b_{p+1}-a_{p+1}-1} {_pF_p}\left( 
\begin{smallmatrix}
a_1, & ..., & a_{p}\\
b_1, & ..., & b_{p}\\
\end{smallmatrix};xt
\right) dt \label{eq: integral_representation_p+1Fp+1_original} \\
\implies {_{p+1}F_{p+1}}\left( 
\begin{smallmatrix}
a_1, & ..., & a_{p+1}\\
b_1, & ..., & b_{p+1}\\
\end{smallmatrix};x
\right) & = & \frac{\Gamma\left(b_{p+1}\right)}{\Gamma\left(a_{p+1}\right) \Gamma\left(b_{p+1}-a_{p+1}\right)} \nonumber \\
& & \hspace{-40pt} \int_0^1 t^{b_{p+1}-a_{p+1}-1} \left(1-t\right)^{a_{p+1}-1} {_pF_p}\left( 
\begin{smallmatrix}
a_1, & ..., & a_{p}\\
b_1, & ..., & b_{p}\\
\end{smallmatrix};x-xt
\right) dt \label{eq: integral_representation_p+1Fp+1_altered}
\end{eqnarray}
if $\text{Re}\left(b_{p+1}\right) > \text{Re}\left(a_{p+1}\right) > 0$. Using the addition formula in Eq. \ref{eq: addition_theorem_pFp},
to replace $_pF_p\left(x-xt\right)$ on the right hand side in Eq. \ref{eq: integral_representation_p+1Fp+1_altered}, given $b_i \neq 0, -1, -2, ...$ $\forall$ $1 \leq i \leq p$, we rewrite it after switching the order of integral and summation as
\begin{eqnarray}
e^x \sum_{j_1=0}^\infty \sum_{j_2=0}^\infty ... \sum_{j_p=0}^\infty \prod_{q=1}^p \frac{\left(b_q - a_q\right)_{j_q} \left(a_q\right)_{u_{q-1}}}{\left(b_q\right)_{u_q}} \frac{\left(-x\right)^{j_q}}{j_q!} \cdot \frac{\Gamma\left(b_{p+1}\right)}{\Gamma\left(a_{p+1}\right) \Gamma\left(b_{p+1}-a_{p+1}\right)} \nonumber \\
\int_0^1 t^{b_{p+1}-a_{p+1}-1} \left(1-t\right)^{a_{p+1}-1} {_{p}}F_{p}\left( 
\begin{smallmatrix}
a_1 + u_0, & a_2 + u_1, & a_3 + u_2, & ..., & a_p + u_{p-1}\\
b_1 + u_1, & b_2 + u_2, & b_3 + u_3, & ..., & b_p + u_p \\
\end{smallmatrix};-xt
\right)dt,
\end{eqnarray}
where $u_0 = 0$, $u_1 = j_1$, $u_2 = j_1+j_2$, ..., $u_p = \sum_{q=1}^p j_q$. Now, using Eq. \ref{eq: integral_representation_p+1Fp+1_original}, the above expression can be rewritten as
\begin{align}
& e^x \sum_{j_1=0}^\infty \sum_{j_2=0}^\infty ... \sum_{j_p=0}^\infty \prod_{q=1}^p \frac{\left(b_q - a_q\right)_{j_q} \left(a_q\right)_{u_{q-1}}}{\left(b_q\right)_{u_q}} \frac{\left(-x\right)^{j_q}}{j_q!} \cdot \nonumber \\
& \hspace{30pt}{_{p+1}}F_{p+1}\left(
\begin{smallmatrix}
a_1 + u_0, & a_2 + u_1, & a_3 + u_2, & ..., & a_p + u_{p-1}, & b_{p+1} - a_{p+1}\\
b_1 + u_1, & b_2 + u_2, & b_3 + u_3, & ..., & b_p + u_p, & b_{p+1} \\
\end{smallmatrix};-x
\right),
\end{align}
which upon substitution in Eq. \ref{eq: integral_representation_p+1Fp+1_altered} yields the Kummer-type transformation in Eq. \ref{eq: kummer_transform_pFp}. While we derived the result requiring that $\text{Re}\left(b_{p+1}\right) > \text{Re}\left(a_{p+1}\right) > 0$, the result holds for general parameters by appealing to analytic continuation, since both sides in Eq. \ref{eq: kummer_transform_pFp} are analytic functions of $a_{p+1}$ and $b_{p+1}$ (given $b_{p+1} \neq 0, -1, -2, ...$). Note that for $p=1$, Eq. \ref{eq: kummer_transform_pFp} is the same as the Kummer transformation derived by Paris (see their Eq. 3) for $_2F_2\left(x\right)$ with general parameters \cite{paris2005kummer}.

\section{Theorem 4 (Euler-type transformation for $\boldsymbol{_{p+1}F_p}$ with general parameters).} For $p \geq 1$ and given $b_i \neq 0, -1, -2, ...$ $\forall$ $1 \leq i \leq p+1$, $\left|x\right| < 1$ and $\text{Re}\left(x\right) < 1/2$,
\begin{align}
& _{p+2}F_{p+1}\left( 
\begin{smallmatrix}
a_0, & a_1 , & a_2 , & a_3 , & ..., & a_p , & a_{p+1}\\
& b_1 , & b_2 , & b_3 , & ..., & b_p , & b_{p+1} \\
\end{smallmatrix};x
\right) \nonumber \\
& = \left(\frac{1}{1-x}\right)^{a_0} \sum_{j_1 = 0}^\infty \sum_{j_2 = 0}^\infty ... \sum_{j_p = 0}^\infty \left(a_0\right)_{u_p} \prod_{q=1}^p \frac{\left(b_q - a_q\right)_{j_q} \left(a_q\right)_{u_{q-1}}}{\left(b_q\right)_{u_q} j_q !} \left(\frac{x}{x-1}\right)^{j_q}  \cdot \nonumber \\
  & \hspace{25pt} {_{p+2}}F_{p+1}\left( 
\begin{smallmatrix}
a_0+u_p, & a_1 + u_0, & a_2 + u_1, & a_3 + u_2, & ..., & a_p + u_{p-1}, & b_{p+1} - a_{p+1}\\
& b_1 + u_1, & b_2 + u_2, & b_3 + u_3, & ..., & b_p + u_p, & b_{p+1} \\
\end{smallmatrix};\frac{x}{x-1}
\right), \label{eq: euler_transform_p+1Fp}
\end{align}
where $u_0 = 0$, $u_1 = j_1$, $u_2 = j_1+j_2$, ..., $u_p = \sum_{q=1}^p j_q$. 

\subsubsection*{Proof}

To prove this theorem, we follow the same approach as that for theorem 3. Hence, we utilize the integral representation for $_{p+2}F_{p+1}\left(x\right)$ function, which is available from Slater \cite{slater1966generalized} (see Eq. 4.8.3.12) and the addition formula for $_{p+1}F_p\left(x\right)$, which we obtained previously (theorem 2). The integral representation for $_{p+2}F_{p+1}\left(x\right)$ is 
\begin{eqnarray}
_{p+2}F_{p+1}\left( 
\begin{smallmatrix}
a_0, & a_1, & ..., & a_{p+1}\\
& b_1, & ..., & b_{p+1}\\
\end{smallmatrix};x
\right) & = & \frac{\Gamma\left(b_{p+1}\right)}{\Gamma\left(a_{p+1}\right) \Gamma\left(b_{p+1}-a_{p+1}\right)} \nonumber \\
& & \hspace{-90pt} \int_0^1 t^{a_{p+1}-1} \left(1-t\right)^{b_{p+1}-a_{p+1}-1} \cdot {_{p+1}F_p}\left( 
\begin{smallmatrix}
a_0, & a_1, & ..., & a_{p}\\
& b_1, & ..., & b_{p}\\
\end{smallmatrix};xt
\right) dt \label{eq: integral_representation_p+2Fp+1_original} \\
\implies {_{p+2}F_{p+1}}\left( 
\begin{smallmatrix}
a_0, & a_1, & ..., & a_{p+1}\\
& b_1, & ..., & b_{p+1}\\
\end{smallmatrix};x
\right) & = & \frac{\Gamma\left(b_{p+1}\right)}{\Gamma\left(a_{p+1}\right) \Gamma\left(b_{p+1}-a_{p+1}\right)}  \nonumber \\
& & \hspace{-90pt} \int_0^1 t^{b_{p+1}-a_{p+1}-1} \left(1-t\right)^{a_{p+1}-1} \cdot {_{p+1}F_p}\left( 
\begin{smallmatrix}
a_0, & a_1, & ..., & a_{p}\\
& b_1, & ..., & b_{p}\\
\end{smallmatrix};x-xt
\right) dt \label{eq: integral_representation_p+2Fp+1_altered}
\end{eqnarray}
if $\text{Re}\left(b_{p+1}\right) > \text{Re}\left(a_{p+1}\right) > 0$ and $\left|x\right| < 1$. Using the addition formula in Eq. \ref{eq: add_theorem_p+1Fp}
on the right hand side in Eq. \ref{eq: integral_representation_p+2Fp+1_altered}, given $b_i \neq 0, -1, -2, ...$ $\forall$ $1 \leq i \leq p$, we rewrite it after switching the order of integral and summation as
\begin{align}
& \left(\frac{1}{1-x}\right)^{a_0} \sum_{j_1 = 0}^\infty \sum_{j_2 = 0}^\infty ... \sum_{j_p = 0}^\infty \left(a_0\right)_{u_p} \prod_{q=1}^p \frac{\left(b_q - a_q\right)_{j_q} \left(a_q\right)_{u_{q-1}}}{\left(b_q\right)_{u_q} j_q !} \left(\frac{x}{x-1}\right)^{j_q} \frac{\Gamma\left(b_{p+1}\right)}{\Gamma\left(a_{p+1}\right) \Gamma\left(b_{p+1}-a_{p+1}\right)} \nonumber \\
& \int_0^1 t^{b_{p+1}-a_{p+1}-1} \left(1-t\right)^{a_{p+1}-1} {_{p+1}}F_{p}\left( 
\begin{smallmatrix}
a_0+u_p, & a_1 + u_0, & a_2 + u_1, & a_3 + u_2, & ..., & a_p + u_{p-1}\\
& b_1 + u_1, & b_2 + u_2, & b_3 + u_3, & ..., & b_p + u_p\\
\end{smallmatrix};\frac{xt}{x-1}
\right)dt.
\end{align}
Now, using the relationship in Eq. \ref{eq: integral_representation_p+2Fp+1_original}, the above expression can be rewritten as
\begin{eqnarray}
\left(\frac{1}{1-x}\right)^{a_0} \sum_{j_1 = 0}^\infty \sum_{j_2 = 0}^\infty ... \sum_{j_p = 0}^\infty \left(a_0\right)_{u_p} \prod_{q=1}^p \frac{\left(b_q - a_q\right)_{j_q} \left(a_q\right)_{u_{q-1}}}{\left(b_q\right)_{u_q} j_q !} \left(\frac{x}{x-1}\right)^{j_q}  \cdot \nonumber \\
 \hspace{25pt} {_{p+2}}F_{p+1}\left( 
\begin{smallmatrix}
a_0+u_p, & a_1 + u_0, & a_2 + u_1, & a_3 + u_2, & ..., & a_p + u_{p-1}, & b_{p+1} - a_{p+1}\\
& b_1 + u_1, & b_2 + u_2, & b_3 + u_3, & ..., & b_p + u_p, & b_{p+1} \\
\end{smallmatrix};\frac{x}{x-1}
\right),
\end{eqnarray}
which upon substitution in Eq. \ref{eq: integral_representation_p+2Fp+1_altered} yields the Euler-type transformation in Eq. \ref{eq: euler_transform_p+1Fp}. Once again, we note by appealing to analytic continuation that the result holds for general parameters even if the requirement $\text{Re}\left(b_{p+1}\right) > \text{Re}\left(a_{p+1}\right) > 0$ is relaxed because both sides in Eq. \ref{eq: euler_transform_p+1Fp} are analytic functions of $a_{p+1}$ and $b_{p+1}$ (given $b_{p+1} \neq 0, -1, -2, ...$).

\section{Summary}

In summary, we have obtained addition formulas for $_pF_p$ and $_{p+1}F_p$ generalized hypergeometric functions, which we utilized to derive Kummer- and Euler-type transformations of these functions with general parameters. The theorems stated herein could be used to derive transformation formulas for other special cases of interest. 

\section*{Declaration of interest}

The author declares that he has no competing interests.

\bibliography{library}

\bibliographystyle{abbrv}

\end{document}